\documentclass[12pt,a4paper]{article}

\usepackage{fancyhdr}
\usepackage{epsfig}
\usepackage{cite}
\usepackage{theorem}
\usepackage{graphicx}
\usepackage{latexsym}
\usepackage{epic}
\usepackage{amsmath}
\usepackage{amsfonts,amssymb,amsopn}
\usepackage[all,2cell,dvips]{xy} \UseAllTwocells \SilentMatrices

\DeclareMathOperator{\degree}{deg}
\DeclareMathOperator{\Ker}{Ker}

\DeclareMathOperator{\Autom}{Aut}

\DeclareMathOperator{\Homom}{Hom}

\DeclareMathOperator{\Rat}{Rat}
\DeclareMathOperator{\sRat}{\underline{\Rat}}
\DeclareMathOperator{\Prin}{Prin}
\DeclareMathOperator{\sPrin}{\underline{\Prin}}
\DeclareMathOperator{\Sym}{Sym}
\DeclareMathOperator{\rank}{rk}
\DeclareMathOperator{\C}{\mathbb{C}}
\DeclareMathOperator{\strs}{{\cal O}_{\textit{X}}}
\DeclareMathOperator{\Spn}{Sp_{\textit{n}}\C}
\DeclareMathOperator{\Spt}{Sp_{2}\C}
\DeclareMathOperator{\Mxn}{{\cal M}_{\textit{X}}(\Spn)}
\DeclareMathOperator{\Mxt}{{\cal M}_{\textit{X}}(\Spt)}

\newcommand{\qed}{\ifhmode\unskip\nobreak\fi\quad\ensuremath\square}

\title{A remark on subbundles of symplectic and orthogonal vector bundles over curves}
\author{George H.\ Hitching}

{\theorembodyfont{\slshape}\newtheorem{prop}{{\textbf Proposition}}}
{\theorembodyfont{\slshape}\newtheorem{thm}[prop]{{\textbf Theorem}}}
{\theorembodyfont{\slshape}}
{\theorembodyfont{\slshape}\newtheorem{cor}[prop]{{\textbf Corollary}}}
{\theorembodyfont{\slshape}\newtheorem{lemma}[prop]{{\textbf Lemma}}}
{\theorembodyfont{\slshape}}
{\theorembodyfont{\slshape}\newtheorem{crit}[prop]{{\textbf Criterion}}}

\begin{document}

\maketitle

\abstract{We review the notions of symplectic and orthogonal vector bundles over curves, and the connection between principal parts and extensions of vector bundles. We give a criterion for a certain extension of rank $2n$ to be symplectic or orthogonal. We then describe almost all of its rank $n$ vector subbundles using graphs of sheaf homomorphisms, and give criteria for the isotropy of these subbundles.}

\section{Introduction}
Throughout this paper, $X$ is a complex projective smooth curve of genus $g$, with structure sheaf $\strs$ and function field $K(X)$.

\subsection{Symplectic and orthogonal vector bundles}
\textbf{Definition:} Let $K$ be a field and $M$ and $P$ vector spaces of dimensions $m$ and $1$ over $K$ respectively. Let $\alpha \colon M \to \Homom(M,P)$ be a $K$--linear map. Then the \textsl{transpose} of $\alpha$ is the unique linear map $^{t}\alpha \colon M \to \Homom(M,P)$ satisfying
\[ \left( \alpha (m_{1}) \right) (m_{2}) = \left( {^{t}\alpha} (m_{2}) \right) (m_{1}) \]
for all $m_{1}, m_{2} \in M$. If $A$ is the matrix of $\alpha$ with respect to some choice of bases, the matrix of $^{t}\alpha$ is just $^{t}A$. The following proposition shows that this also makes sense for vector bundles.
\begin{prop}
Let $E \to X$ be a vector bundle, $L \to X$ a line bundle and $\alpha \colon E \to \Homom(E,L)$ a vector bundle map over an open set $U \subseteq X$. Then the transpose of $\alpha$ is a well--defined map $E \to \Homom(E,L)$.
\label{gensym}
\end{prop}
\textbf{Proof}\\
Suppose the bundles $E$ and $L$ have transition functions $\{ e_{i,j} \}$ and $\{ l_{i,j} \}$ respectively relative to an open cover $\{ U_{i}: i \in J \}$ of $X$. Then $\alpha$ is given by a cochain $\{ \alpha_{i} \}$ of $n \times n$ matrices which satisfy
\begin{equation} \alpha_{i} e_{i,j} = \left( ^{t}e_{i,j}^{-1} l_{i,j} \right) \alpha_{j} \label{homcond} \end{equation}
over $U \cap U_{i} \cap U_{j}$. If $^{t}\alpha$ exists then it should be given by the cocycle $\{ {^{t}\alpha_i} \}$, by the discussion before this proposition. Taking the transpose of (\ref{homcond}), we have
\[ {^{t}e_{i,j}} \left( ^{t}\alpha_{i} \right) = {^{t}\alpha_j} l_{i,j} e_{i,j}^{-1}, \]
equivalently
\[ {^{t}\alpha_i} e_{i,j} = \left( ^{t}e_{i,j}^{-1} l_{i,j} \right) {^{t}\alpha_j} \]
since the $l_{i,j}$ commute with the other transition functions. This shows that $\{ {^{t}\alpha_i} \}$ does indeed define a map $E \to \Homom(E,L)$ over $U$. \qed \\
\\
By this proposition, it makes sense to speak of symmetric and antisymmetric homomorphisms $E \to \Homom(E,L)$. We denote these by $\Sym(E, \Homom(E,L))$ and $\bigwedge(E, \Homom(E,L))$ respectively.\\
\\
Similar statements hold for maps $\Homom(E,L) \to E$.\\
\\
\textbf{Definition:} A vector bundle $W \to X$ is \textsl{symplectic} (resp., \textsl{orthogonal}) if there exists a bilinear nondegenerate antisymmetric (resp., symmetric) form $\theta$ on $W \times W$ with values in a line bundle $L$.\\
\\
Two immediate consequences of the nondegeneracy of $\theta$ are:
\begin{itemize}
\item There is an antisymmetric or symmetric isomorphism
\[ W \xrightarrow{\sim} \Homom(W, L) \]
given by $w \mapsto \theta(w, \cdot)$. In particular, $(\det W)^{2} = L^{\rank W}$.
\item $W$ is symplectic only if it has even rank, since skew--symmetric matrices have even rank.
\end{itemize}
We shall henceforth restrict ourselves to the case where $\rank W = 2n$, even if $W$ is orthogonal.
\par
A subbundle of $W$ is \textsl{isotropic} if $\theta$ restricts to zero on it. For any subbundle $E \subseteq W$, we have the short exact vector bundle sequence
\[ 0 \to E^{\perp} \to W \to \Homom(E, L) \to 0  \label{orthcomp} \]
where the surjection is the map $w \mapsto \theta( w , \cdot )|_E$ and
\[ E^{\perp} = \left\{ w \in W : \theta( w , E ) = 0 \right\} \]
is the \textsl{orthogonal complement} of $E$ with respect to $\theta$. Clearly $E$ is isotropic if and only if $E \subseteq E^{\perp}$; this shows that the rank of an isotropic subbundle is at most $n = \frac{1}{2} \rank W$. An isotropic subbundle of rank $n$ is called a \textsl{Lagrangian} subbundle.

\subsection{Principal parts and extensions}

For vector bundles $E$ and $F$ over $X$, it is well known that an extension
\begin{equation} 0 \to E \to W \to F \to 0 \label{FextE} \end{equation}
is determined up to isomorphism of extensions by its cohomology class $\delta(W) \in H^{1}(X, \Homom(F,E))$. This class can be realised explicitly as follows. There exist trivialisations of $W$ over some open cover $\{ U_{i} \}$ of $X$ whose transition function over the intersection $U_{i} \cap U_j$ is of the form
\[ \begin{pmatrix} e_{i,j} & \delta_{i,j} \\ 0 & f_{i,j} \end{pmatrix} \]
where $e_{i,j}$ and $f_{i,j}$ are transition functions for $E$ and $F$ respectively. Then $\delta(W)$ is the class determined by the $\Homom(F,E)$--valued $1$--cocycle $\{ \delta_{i,j} \}$.
\par
We give another description of $\delta(W)$ following Kempf \cite{Kem1983}. (Here the results are given for extensions of invertible sheaves, but the arguments are readily adapted to the case of arbitrary rank.) Firstly, we fix some notation. We denote the sheaf of regular sections of a vector bundle $W$, $E$, $L$ etc.\ by the corresponding script letter ${\cal W}$, ${\cal E}$, ${\cal L}$ etc. A vector bundle $E \to X$ gives rise to an exact sequence of $\strs$--modules
\[ 0 \to {\cal E} \to \sRat(E) \to \sPrin(E) \to 0 \]
where $\sRat(E)$ is the sheaf of rational sections of $E$ and $\sPrin(E)$ the sheaf of principal parts with values in $E$. We denote their groups of global sections by $\Rat(E)$ and $\Prin(E)$ respectively. The sheaves $\sRat(E)$ and $\sPrin(E)$ are flasque, so we have the cohomology sequence
\begin{equation} 0 \to H^{0}(X, E) \to \Rat(E) \to \Prin(E) \to H^{1}(X, E) \to 0. \label{cohomseq} \end{equation}
We denote $\overline{s}$ the principal part of $s \in \Rat(E)$, and we write $[p]$ for the class in $H^{1}(X, E)$ of $p \in \Prin(E)$.
\par
Now given an extension of the form (\ref{FextE}), it can be shown (see for example \cite{Hit2005}, Lemma 3.1) that there exists a unique principal part $p \in \Prin(\Homom(F,E))$ such that the sheaf of sections ${\cal W}$ is of the form
\[ {\cal W}_{p} := \left\{ (e,f) \in \sRat(E) \oplus \sRat(F) : f \hbox{ is regular and } \overline{e} = p(f) \right\}. \]
Following Kempf \cite{Kem1983}, Chap.\ 6, one proves as for the case where ${\cal E}$ and ${\cal F}$ are invertible that two principal parts define isomorphic extensions if and only if they differ by the principal part $\overline{\alpha}$ of some $\alpha \in \Rat(\Homom(F,E))$, and that the cohomology class $\delta(W)$ is just $[p]$.\\
\par
In the following sections, we give a criterion for a certain extension of vector bundles to be symplectic or orthogonal. We then generalise a result from Mukai \cite{Muk2001} to describe almost all rank $n$ subbundles of such an extension, and give criteria for the isotropy of these subbundles. We conclude by sketching how these results are applicable to the study of moduli spaces of symplectic or orthogonal vector bundles over curves.\\
\par
\textbf{Acknowledgements:} I thank my doctoral supervisors C.\ Pauly, J.\ Bolton, W.\ Klingenberg and W.\ Oxbury for their time and ideas. I am grateful to S.\ Ramanan for showing me Criterion \ref{fundam}, which is central to this work. Thanks to M.\ Reid for his continuing support, mathematical and practical. I acknowledge gratefully the financial support and hospitality of the University of Durham and l'Universit\'e de Nice et Sophia--Antipolis.

\section{Symplectic and orthogonal extensions}

Let $W \to X$ be a symplectic or orthogonal vector bundle of rank $2n$ and let $E \subset W$ be a Lagrangian subbundle. Then $W$ is an extension
\[ 0 \to E \to W \to \Homom(E, L) \to 0. \tag{$\delta(W)$} \]

Conversely, it is natural to ask for which extension classes $\delta(W)$ this sequence is induced by a bilinear antisymmetric or symmetric form. We have

\begin{crit} An extension $0 \to E \to W \to \Homom(E, L) \to 0$ has a symplectic (resp., orthogonal) structure with respect to which $E$ is isotropic if and only if $W$ is isomorphic as a vector bundle to an extension whose cohomology class belongs to $H^{1}(X, \Sym(\Homom(E,L),E))$ (resp., $H^{1}(X, \bigwedge(\Homom(E,L), E))$). \label{fundam}
\end{crit}
\textbf{Proof}\\
We prove the criterion for the symplectic case; the orthogonal case is practically identical.
\par
$\Leftarrow$: Firstly, suppose $\delta(W)$ is actually symmetric. By the discussion in $\S$ 1.2, there exists $p \in \Prin(\Homom(\Homom(E, L), E))$ such that the sheaf ${\cal W}$ is equal to
\begin{equation} \left\{ (f,\phi) \in  \sRat(E) \oplus \sRat(\Homom(E,L)) : \phi \hbox{ is regular and } \overline{f} = p(\phi) \right\}. \label{sectW} \end{equation}
To say that $\delta(W) = [p]$ is symmetric is to say that
\[ ^{t}p - p = \overline{\alpha} \]
for some $\alpha \in \Rat( \Homom( \Homom(E,L), E))$. Clearly $\overline{\alpha}$ is antisymmetric; replacing $\alpha$ by $\frac{\alpha - {^{t}\alpha}}{2}$ if necessary, we can assume that $\alpha$ itself is antisymmetric.
\par
Now we define a $\sRat(L)$--valued bilinear antisymmetric form
\[ \underline{\theta} \colon \left( \sRat(E) \oplus \sRat(\Homom(E,L)) \right)^{\times 2} \to \sRat(L) \]
by setting
\[ \underline{\theta} \left( (e_{1},\phi_{1}) , (e_{2},\phi_{2}) \right) = \phi_{1}(e_{2}) - \phi_{2}(e_{1}) - \phi_{2} \left( \alpha (\phi_{1}) \right). \]
By the description in (\ref{sectW}), for any $(e_{1},\phi_{1})$, $(e_{2},\phi_{2}) \in {\cal W}_p$, the principal part
\begin{align*} \overline{\underline{\theta} \left( (e_{1},\phi_{1}) , (e_{2},\phi_{2}) \right) } &= \phi_{1}(p(\phi_{2})) - \phi_{2}(p(\phi_{1})) - \phi_{2} \left( \overline{\alpha} ( \phi_{1}) \right) \\
&= \phi_{2}((^{t}p - p - \overline{\alpha})(\phi_{1}))\\
&= 0
\end{align*}
since $^{t}p - p = \overline{\alpha}$. Hence $\underline{\theta}$ is regular on ${\cal W} \times {\cal W}$. It is clearly $\strs$--bilinear and nondegenerate, and ${\cal E}$ is isotropic. Thus $\underline{\theta}$ induces a global regular symplectic form $\theta \colon W \times W \to L$ with the required properties.
\par
For the general case, we note that this form pulls back to give the required symplectic structure to any vector bundle isomorphic to $W$, which need not be isomorphic as an extension.\\
\par
$\Rightarrow$: Choose transition functions $\{ e_{i,j} \}$ and $\{ l_{i,j} \}$ for $E$ and $L$ respectively over an open cover $\{ U_{i}: i \in J \}$ of $X$. Then the transition functions of $\Homom(E, L)$ are $\{ {^{t}e_{i,j}^{-1}}l_{i,j} \}$ and there exist trivialisations for $W$ over $\{ U_{i} \}$ whose transition functions are of the form
\[ \{ w_{i,j} \} = \left\{ \begin{pmatrix} e_{i,j} & \delta_{i,j} \\ 0 & {^{t}e_{i,j}^{-1}} l_{i,j} \end{pmatrix} \right\}; \]
the cohomology class of the cocycle $\{\delta_{i,j}\}$ is $\delta(W)$.
\par
The symplectic form is given with respect to $\{ U_{i} \}$ by a cochain $\{\Theta_{i}\}$ of antisymmetric matrices which satisfy
\begin{equation} (l_{i,j}^{-1}){^{t}w_{i,j}}\Theta_{i}w_{i,j} =  \Theta_{j} \label{coch} \end{equation}
on the intersection $U_{i} \cap U_{j}$ for all $i, j \in J$, since the symplectic form defines a homomorphism $W \to \Homom(W, L)$. We can write
\[ \Theta_{i} = \begin{pmatrix} A_i & B_i \\ -{^{t}B_i} & C_i \end{pmatrix} \]
where $\{A_{i}\}$, $\{B_{i} \}$ and $\{C_i\}$ are $M_{n,n}(\C)$--valued cochains and all the $A_i$ and $C_i$ are antisymmetric. Firstly, we see that every $A_i \equiv 0$ because $E \subset W$ is isotropic.
\par
Expanding condition (\ref{coch}), we see that
\[ B_{i} \left( ^{t}e_{i,j}^{-1} \right) = \left( ^{t}e_{i,j}^{-1} \right) B_{j}, \]
so $\{ B_{i} \}$ defines an endomorphism of $E^{*}$ and also of $\Homom(E, L)$. Since all the $A_{i}$ are zero but the form is nondegenerate, this must be an automorphism. Also by (\ref{coch}), we have
\begin{equation} {^{t}\delta_{i,j}} B_{i} \left( ^{t}e_{i,j}^{-1} \right) - e_{i,j}^{-1} \left(^{t}B_{i} \right) \delta_{i,j} = C_{j} - e_{i,j}^{-1}C_{i} \left( ^{t}e_{i,j}^{-1} l_{i,j} \right), \label{transcomp} \end{equation}
whence we see that the difference between the cocycle $\{ {^{t}B_i} \delta_{i,j} \}$ and its transpose is cohomologically trivial. Hence the cohomology class defined by $\{ {^{t}B_i}\delta_{i,j} \}$ belongs to $H^{1}(X, \Sym(\Homom(E,L),E))$. This belongs to the same orbit as $\delta(W)$ under the action of $\Autom E$ on the extension space $H^{1} \left( X, \Homom(\Homom(E,L),E) \right)$, so defines an extension isomorphic to $W$ as a vector bundle. \qed\\
\\
\textbf{Caution:} Strictly, in order to obtain cocycles which map between the correct spaces in (\ref{transcomp}), we should replace ${^{t}\delta_{i,j}} B_{i} (^{t}e_{i,j}^{-1})$ with
\[ l_{i,j}^{-1}(^{t}\delta_{i,j}) B_{i} (^{t}e_{i,j}^{-1} l_{i,j}), \]
but this does not change the value. In any case, by Prop.\ \ref{gensym} the transpose of $\{ e_{i,j}^{-1} (^{t}B_{i}) \delta_{i,j} \}$ is indeed defined by the transposed cocycle $\{ {^{t}\delta_{i,j}} B_{i} (^{t}e_{i,j}^{-1}) \}$.\\
\\
\textbf{Remark:} If $E$ is simple (for example, if $E$ is stable) then the cocycle $\{ \delta_{i,j} \}$ will itself be symmetric because $\{ B_{i} \}$ defines a nonzero homothety.

\newpage

\section{Vector subbundles and graphs}
Firstly, we recall some linear algebra. Let $K$ be a field and suppose $M$, $N$ and $P$ are vector spaces over $K$ of dimensions $m$, $n$ and $1$ respectively. If $N = \Homom(M,P)$ and
\[ \alpha \colon \Homom(M,P) \to M \]
is an antisymmetric map then then we can define a $P$--valued bilinear nondegenerate antisymmetric form
\[ \theta \colon \left( M \oplus \Homom(M, P) \right)^{\times 2} \to P \]
by
\[ \theta \left( (m_{1},\psi_{1} ),(m_{2},\psi_{2} ) \right) = \psi_{1}(m_{2}) - \psi_{2}(m_{1}) - \psi_{2} \left( \alpha ( \psi_{1} ) \right) \]
The following is a slight generalisation of Mukai \cite{Muk2001}, Example 1.5.
\begin{lemma} \label{ano}
\begin{enumerate}
\renewcommand{\labelenumi}{$\mathrm{(\roman{enumi})}$}
\item There is a bijection between $\Homom_{K}( N , M)$ and the set of $n$--dimensional $K$--vector subspaces of $M \oplus N$ intersecting $M$ in zero, given by associating to a map $\beta$ its graph $\Gamma_{\beta}$.
\item The kernel of $\beta$ is canonically isomorphic to $\Gamma_{\beta} \cap \left( \{ 0 \} \oplus N \right)$.
\item If $N = \Homom(M,P)$ then $\Gamma_{\beta}$ is isotropic with respect to $\theta$ if and only ${^{t}\beta} - \beta = \alpha$.
\end{enumerate}
\end{lemma}
\textbf{Proof}\\
This is straightforward to check. \qed \\
\\
Now let $E$ and $F$ be vector bundles of rank $m$ and $n$ respectively over $X$. Consider an extension $0 \to E \to W \to F \to 0$ with sheaf of sections ${\cal W}_p$ and class $\delta(W) = [p] \in H^{1}(X, \Homom(F,E))$. We want to study vector subbundles $G \subset W$ of rank $n$ whose projection to $F$ is generically surjective.\\
\\
\textbf{Notation:} Let $\Lambda$ be a vector space over $K(X)$. By analogy with Hartshorne \cite{Har1977}, p.\ 69, we write $\underline{\Lambda}$ for the constant sheaf on $X$ associated to $\Lambda$.\\
\\
Now $\Rat(E)$ and $\Rat(F)$ are vector spaces of dimensions $m$ and $n$ respectively over $K(X)$, the field of rational functions on $X$. The following theorem, globalising Lemma \ref{ano}, is a generalisation of Mukai \cite{Muk2001}, Example 1.7, to the case where $W$ may be a nontrivial extension.

\newpage

\begin{thm}
\begin{enumerate}
\renewcommand{\labelenumi}{$\mathrm{(\roman{enumi})}$}
\item There is a bijection
\[ \Homom_{K(X)}(\Rat(F),\Rat(E)) \leftrightarrow \left\{ \begin{array}{c} \hbox{rank $n$ vector subbundles} \\
G \subset W \hbox{ with } G|_{x} \cap E|_{x} = 0 \\
\hbox{ for generic } x \in X \end{array} \right\}  \]
given by $\beta \leftrightarrow \underline{\Gamma_{\beta}} \cap {\cal W}_{p} =: {\cal G}_{\beta}$. Moreover, $\cal{G}_{\beta}$ is isomorphic to the kernel of the map
\[ \left( p - \overline{\beta} \right) \colon {\cal F} \to \sPrin(E). \]
\item Let $\widetilde{\beta}$ denote the restriction of $\beta$ to $\Ker ( p - \overline{\beta} ) \subseteq {\cal F}$. Then $\Ker ( \widetilde{\beta} ) \cong ( \underline{\Gamma_{\beta}} \cap {\cal W}_{p} ) \cap (\{ 0 \} \oplus {\cal F})$ (although we will not use this result).
\item Suppose $F = \Homom(E,L)$ and ${^{t}p} - p = \overline{\alpha}$ for some
\[ \alpha \in \Rat \left( \bigwedge \left( \Homom(E,L), E \right) \right), \]
so $W$ carries the symplectic form $\theta$ defined in the proof of Criterion \ref{fundam}. Then $G_{\beta}$ is isotropic with respect to $\theta$ if and only if ${^{t}\beta} - \beta = \alpha$. \label{vbgraphs}
\end{enumerate}
\end{thm}
\textbf{Proof}\\
(i) Let $G \subset W$ be a vector subbundle of rank $n$ intersecting $E$ in zero except at a finite number of points. Then $\Rat(G)$ is a $K(X)$--vector subspace of $\Rat(W) = \Rat(E) \oplus \Rat(F)$ of dimension $m$ which intersects $\Rat(E)$ in zero. By Lemma \ref{ano} and the remarks just before this theorem, $\Rat(G)$ is the graph $\Gamma_{\beta}$ of some uniquely determined
\[ \beta \in \Homom_{K(X)} \left( \Rat(F), \Rat(E) \right).\]
Furthermore, ${\cal G} = \underline{\Gamma_{\beta}} \cap {\cal W}_p$ since a regular section of $G$ is the same thing as a rational section of $G$ which is a regular section of $W$.
\par
Conversely, we claim that the the association
\[ \beta \mapsto \underline{\Gamma_{\beta}} \cap {\cal W}_{p} =: {\cal G}_{\beta} \]
defines a subsheaf of ${\cal W}_p$ which in fact corresponds to a vector subbundle $G_{\beta} \subset W$ with the required properties. By the definitions of $\Gamma_{\beta}$ and ${\cal W}_p$, we have
\[ {\cal G}_{\beta} = \left\{ ( \beta(f) , f ) \in \sRat(E) \oplus {\cal F} : p(f) = \overline{\beta(f)} \right\}. \]
This is clearly isomorphic to the kernel of the map
\[ (p - \overline{\beta}) \colon {\cal F} \to \sPrin(E) \]
via the projection of ${\cal G}_{\beta}$ onto its image in ${\cal F}$. (The inverse map is $f \mapsto (\beta(f), f)$.) But since any principal part is supported at a finite number of points, $({\cal G}_{\beta})_{x} \cong {\cal F}_x$ for all but finitely many $x \in X$. Hence ${\cal G}_{\beta}$ has rank $n$ and projects surjectively to ${\cal F}$ at all but a finite number of points of $X$.
\par
We now check that the inclusion ${\cal G}_{\beta} \hookrightarrow {\cal W}_p$ actually corresponds to a vector bundle injection $G_{\beta} \hookrightarrow W$. We have a short exact sequence of $\strs$--modules
\[ 0 \to {\cal G}_{\beta} \to {\cal W}_{p} \to {\cal Q} \to 0 \]
where ${\cal Q}$ is coherent. Let ${\cal G}^{\prime}$ denote the inverse image in ${\cal W}_p$ of the torsion subsheaf of ${\cal Q}$. Clearly ${\cal G}^{\prime}$ contains ${\cal G}_{\beta}$. Now ${\cal G}^{\prime}$ corresponds to an injection of vector bundles $G^{\prime} \hookrightarrow W$ by Atiyah \cite{Ati1957}, Prop.\ 1, since by construction ${\cal W}_{p}/{\cal G}^{\prime}$ is locally free. But in fact $\Rat(G^{\prime}) = \Gamma_{\beta}$; this is because ${\cal G}^{\prime}$ is contained in ${\cal G}_{\beta}(D)$ for some divisor $D$ on $X$, so they have the same sheaf of rational sections. Hence ${\cal G}^{\prime} \subseteq \Gamma_{\beta} \cap {\cal W}_{p} = {\cal G}_{\beta}$, so ${\cal G}_{\beta} = {\cal G}^{\prime}$ corresponds to a vector subbundle $G_{\beta} \subset W$.
\par
It is not hard to check that these constructions are mutually inverse, so we have a bijection.\\
\\
(ii) Suppose $g \in \Ker (p - \overline{\beta} )$. Then $\widetilde{\beta}(g) = 0$ if and only if
\[ (\beta(g), g) \in (\underline{\Gamma_{\beta}} \cap {\cal W}_{p}) \cap ( \{ 0 \} \oplus {\cal F}). \]
 \\
(iii) The symplectic form on $W$ is induced by the restriction of the form defined earlier
\[ \underline{\theta} \colon \sRat \left( E \oplus \Homom(E, L) \right)^{\times 2} \to \sRat(L) \]
to ${\cal W}_{p} \times {\cal W}_p$, so the criterion for isotropy follows from part (iii) of Lemma \ref{ano}. \qed \\
\\
We make an observation:
\begin{lemma}
The $K(X)$--linear map $\beta$ is everywhere regular on ${\cal G}_{\beta}$. \label{reg}
\end{lemma}
\textbf{Proof}\\
If the supports of $p$ and $\overline{\beta}$ are disjoint, then this is clear. Suppose the supports coincide at a point $x \in X$. Then the maps
\[ p_{x}, \overline{\beta}_x \hbox{ and } (p - \overline{\beta})_{x} \in \Homom_{\strs_{,x}}(({\cal G}_{\beta})_{x},\sPrin(E_{\beta})_{x}) \]
are given locally by matrices of rational functions on a neighbourhood of $x$. Since $X$ is of dimension $1$, we can assume that the numerators and denominators of each of these functions are relatively prime, and then in fact the denominators determine the maps.
\par
The key point is that by the identity
\[ \frac{a}{f} - \frac{b}{h} = \frac{ah-bf}{fh}, \]
the denominators of the entries of the matrix $(p - \overline{\beta})_x$ are at worst the products of the corresponding entries of $p_x$ and $\overline{\beta}_x$. Since ${\cal G}_{\beta} \cong \Ker(p - \overline{\beta})$, the value $(p - \overline{\beta})_{x}(g)$ is regular for any $g \in ({\cal G}_{\beta})_x$. But for regular functions $a$, $f$ and $h$, if $\frac{a}{fh}$ is regular then so is $\frac{a}{h}$. Hence $\beta$ itself is regular on $({\cal G}_{\beta})_x$. \qed\\
\\
Now we examine a useful special case.

\begin{cor} Suppose that $h^{0}(X, \Homom(F,E)) = 0$. Then principal parts defining the cohomology class $\delta(W) = [p]$ are in bijection with elementary transformations of $F$ lifting to rank $n$ subbundles of $W$ via $q \leftrightarrow \Ker \left( q \colon {\cal F} \to \sPrin(E) \right)$. \label{bijection} \end{cor}
\textbf{Proof}\\
If $\delta(W)$ is defined by $q \in \Prin ( \Homom(F,E) )$ then by (\ref{cohomseq}) we have $q = p - \overline{\beta}$ for some global rational section $\beta$ of $\Homom(F,E)$, which is uniquely determined by hypothesis. Then $\Ker ( q ) \subseteq {\cal F}$ lifts to the rank $n$ subsheaf $\Gamma_{\beta} \cap {\cal W}_p$ of ${\cal W}_p$ by the map $f \mapsto (\beta(f), f)$. By the proof of Thm.\ 4 (i) this corresponds to a vector subbundle.
\par
Conversely, suppose $G \subset W$ is a rank $n$ vector subbundle lifting from $F$. By the proof of Thm.\ 4 (i), the sheaf of sections of $G$ is isomorphic to $\Ker ( p - \overline{\beta} )$ for some $\beta \in \Rat(\Homom(F,E))$. But by (\ref{cohomseq}) we have $[ p - \overline{\beta} ] = [p] = \delta(W)$.
\par
It is easy to see that these constructions are mutually inverse. \qed

\subsection*{Another criterion for isotropy}
We give a refinement of Theorem 4 (iii) for this case. Let $E$ and $L$ be vector bundles of ranks $n$ and $1$ respectively over $X$ such that $\Homom( \Homom(E, L) , E)$ has no global sections. Let
\[ 0 \to E \to W \to \Homom(E,L) \to 0 \]
be a symplectic extension with sheaf of sections ${\cal W}_p$ where as before ${^{t}p} - p = \overline{\alpha}$. Let $G \subset W$ be a subbundle of rank $n$ which intersects $E$ generically in rank $0$; by Cor.\ \ref{bijection}, the sheaf ${\cal G} = \Ker(q)$ for some $q \in \Prin(\Homom(\Homom(E,L),E))$ such that $\delta(W) = [q]$.

\begin{crit} The subbundle $G \subset W$ is isotropic if and only if $q$ is a symmetric principal part (note that this is a stronger condition than that the cohomology class defined by $q$ be symmetric). \label{isotcrit}
\end{crit}
\textbf{Proof}\\
By Cor.\ \ref{bijection} (i) we have ${\cal G} = \underline{\Gamma_{\beta}} \cap {\cal W}_p$ for some
\[ \beta \in \Rat(\Homom( \Homom(E,L),E)) \]
which, by hypothesis, is determined by the condition $q = p - \overline{\beta}$.
\par
We claim that ${^{t}\beta} - \beta = \alpha$ if and only if $\overline{{^{t}\beta} - \beta} = \overline{\alpha}$. One direction is clear. Conversely, suppose $\overline{{^{t}\beta} - \beta} = \overline{\alpha}$. By (\ref{cohomseq}), the difference ${^{t}\beta} - \beta - \alpha$ is a global regular section of $\Homom( \Homom (E, L) , E)$. But this is zero by hypothesis, so ${^{t}\beta} - \beta = \alpha$.
\par
Thus by Theorem 4 (iii), the subbundle $G$ is isotropic if and only if $\overline{{^{t}\beta} - \beta} = \overline{\alpha}$. Now ${^{t}p} - p = \overline{\alpha}$, so
\[ {^{t}q} - q = {^{t}(p - \overline{\beta})} - (p - \overline{\beta}) = \overline{\alpha} - \left( \overline{^{t}\beta} - \overline{\beta} \right). \]
Hence $G$ is isotropic if and only if $q$ is a symmetric principal part. \qed

\section*{The orthogonal case}
There are obvious analogues to these results for the orthogonal case which are proven identically. In Lemma \ref{ano} (iii), we consider a symmetric map $\alpha \colon \Homom(E,L) \to E$ and work with the standard bilinear nondegenerate \textsl{symmetric} form $\theta$ given by
\[ \theta \left( (v_{1},\phi_{1}),(v_{2},\phi_{2}) \right) = \phi_{1}(v_{2}) + \phi_{2}(v_{1}) - \phi_{2} \left( \alpha ( \phi_{1} ) \right) \]
and require ${^{t}\beta} + \beta = \alpha$. In Theorem 4 (iii), we consider an orthogonal extension ${\cal W}_p$ where ${^{t}p} + p = \overline{\alpha}$ and the condition of the criterion is that ${^{t}\beta} + \beta = \alpha$. With this setup, we get an orthogonal version of Criterion \ref{isotcrit}: the sheaf $G$ is isotropic if and only if $q$ is an antisymmetric principal part.

\section{Applications}

\subsection{Bundles of rank 2}
As an example, we show how these results apply to the well--known case where $n=1$.
\par
Firstly, let $W \to X$ be of rank $2$ and trivial determinant. Then any line subbundle $E \subset W$ gives a short exact sequence
\[ 0 \to E \to W \to E^{*} \to 0 \]
and Criterion \ref{fundam} gives another proof of the well--known fact that every such $W$ is symplectic, since in this case $\Sym^{2}E = \Homom(E^{*},E)$. Furthermore, Theorem 4 (iii) then gives another proof that every line subbundle of $W$ is isotropic.
\par
On the other hand, suppose $W$ is an orthogonal vector bundle of rank $2$ and let $E \subset W$ be an isotropic line subbundle. Then $W \cong E \oplus E^{-1}L$ by Criterion \ref{fundam} since $\bigwedge^{2}E = 0$.

\subsection{Covers of moduli spaces}

Let $\Mxn$ denote the moduli space of semistable principal $\Spn$--bundles over $X$. This is an irreducible projective variety of dimension $n(2n+1)(g-1)$; see Ramanathan \cite{Rthn1996a}, \cite{Rthn1996b} for results on moduli of principal $G$--bundles over curves. $\Mxn$ is naturally a moduli space for semistable vector bundles of rank $2n$ which carry an $O_X$--valued symplectic form (see \cite{Hit2005} for details). Following Narasimhan and Ramanan \cite{NR1984}, one might try to construct a generically finite cover of $\Mxn$ by the classifying map from the union of the extension spaces
\[ \mathbb{P} H^{1}(X, \Sym^{2}E) \]
as $E$ varies over some collection of rank $n$ vector bundles over $X$. By Criterion \ref{fundam}, determining the fibre over a bundle $W$ in the image of such a classifying map involves asking for certain Lagrangian subbundles of $W$. If $n \geq 2$ then not every rank $n$ subbundle need be Lagrangian, and Criterion \ref{isotcrit} may be of use in distinguishing those that are.
\par
In a forthcoming paper (and see \cite{Hit2005}, Chaps.\ 4 \& 5), we shall use this technique to construct a cover of $\Mxt$ when $X$ has genus $2$.\\
\\
\textbf{Remark:} The hypothesis $h^{0}(X,\Homom(F, E)) = 0$ of Cor.\ \ref{bijection} is a natural one if we want to build stable bundles. Firstly, if $\degree E \geq \degree F$ then no extension
\[ 0 \to E \to W \to F \to 0 \]
can be stable. If the converse is true and furthermore $E$ and $F$ are themselves stable then indeed $h^{0}(X,\Homom(F, E)) = 0$.

Laboratoire J.--A.\ Dieudonn\'e,\\
Universit\'e de Nice et Sophia--Antipolis,\\
Parc Valrose,\\
06108 Nice CEDEX 02,\\
France.\\
E--mail: \texttt{hitching@math.unice.fr}

\end{document}